\documentclass{amsart}

\begin{document}

\title{Graded Lie algebras defined by Jordan algebras
and their representations}
\thanks{This work was supported by the RFBR grant 03--01--00056
and Swedish Academy of Sciences.}
\author{Issai Kantor}
\address{Department of Mathematics, Lund University, 
S-221 00 Sweden}
\email{E-mail: kantor@maths.lth.se} 
\author{Gregory Shpiz}
\address{Centre for Continuous Mathematical Education, Moscow}
\email{E-mail: shpiz@theory.sinp.msu.ru}
\date{}

\maketitle

\begin{abstract}

In this talk we introduce the notion of a generalized
representation of a
  Jordan algebra with unit which has the following properties:

  1) Usual representations and Jacobson representations correspond to special cases 
of generalized representations.

2) Every simple Jordan algebra has infinitely many nonequivalent
 generalized representations.

3) There is a one-to-one correspondence between irreducible generalized representations 
of a Jordan algebra $A$ and irreducible representations of a graded Lie algebra 
$L(A)=U_{-1}\oplus U_0\oplus U_1$ 
corresponding to $A$ (the Lie algebra $L(A)$  coincides with
the TKK construction
when $A$ has a 
unit).

The latter correspondence allows to use the theory of representations of Lie algebras to study 
generalized representations of Jordan algebras. In particular, one can classify irreducible 
generalized representations of  semisimple
Jordan algebras and also  obtain classical results about usual
representations and Jacobson representations in a simple way.
\end{abstract}

\vskip 3mm

\begin{center}
\bf Introduction
\end{center}
\par

Jordan algebras were introduced by P.~Jordan, J.~von Neumann and E.~Wigner
 (see [1]) in the
 connection with some problems of quantum mechanics.  Already there  it
 was  found that 
 simple Jordan algebras have only a finite number of nonequivalent
 irreducible representations (homomorphisms
  into a space of linear operators with operation $X*Y=XY+YX$).
Moreover (A.~Albert [2]), the exceptional Jordan algebra $E_{3}$ has no
such representations at all.

This situation demonstrates  a big difference between Lie and Jordan algebras.
(As is known a Lie algebra has infinitely many nonequivalent irreducible
 representations.)

 To improve the situation, N.~Jacobson introduced [3] another notion of a
 representation of a Jordan algebra (see below the definition of the Jacobson representation) 
 and showed that every simple  Jordan algebra has at least one nontrivial Jacobson 
 representation.  But still the number of Jacobson representations is also finite. 

In this talk we will introduce a notion of a generalized representation of
 a Jordan algebra with unit and will describe the irreducible generalized
 representations.  The
 irreducible generalized representations are in a one-to-one
 correspondence with the irreducible
 representations  of the  3-graded Lie algebra $L(A)$ 
corresponding to $A$.

In particular, every simple Jordan algebra has infinitely many nonequivalent irreducible
 generalized representations.

  The usual   and Jacobson representations correspond to special cases of generalized 
representations.  Moreover, this correspondence preserves the irreducibility and the 
equivalence of representations.  In particular, it allows to classify irreducible usual and 
Jacobson representations.

 The authors are very grateful to Bruce Allison, Kevin McCrimmon and Ivan Shestakov
 for useful discussions.

\vskip 3mm
\begin{center}
{\bf \S 1\ A graded Lie algebra $L(A)$ defined by a Jordan algebra $A$}
\end{center}

We need a construction of a 3-graded Lie algebra 
$L(A)$ defined by a Jordan algebra $A$.

 The construction of  $L(A)$ is
 presented as it was originally given in [4], [5] (see also [6]).

This construction coincides with what is called 
the TKK construction when Jordan algebras $A$ has a unit, but does not
coincide with it in general (for example
 $dim U_{-1}$ is not equal in general to $ dim U_{1}$ ).  The Lie
 algebra  $L(A)$ has the following two  important properties:

 1) There is an element $\bar {A} \in  U_{1} $ such that
 $[[\bar {A},x],y]=x*y\;\;\forall x,y \in   U_{-1} $, where $*$ is 
the multiplication in the given algebra $A$.
 (The space  $ U_{-1} $ is identified with the space of the algebra $A$.)

2) The Lie algebra  $L(A)$ is generated by the space $ U_{-1} $ and
the element  $\bar {A}\in  U_{1} $.

To construct the Lie algebra 
$L(A)$ let us denote by ${\mathcal {U}} $ the space of the algebra $A$ and by $\bar {A}(x,y)=x*y $ 
the multiplication in $A$. We denote also 
$$L_a(x)=a*x,\; \;\;A_a(x,y)=(x*a)*y+(y*a)*x-a*(x*y),\;(0.1)$$

$${\mathcal {S}}=\{L_a,[L_a,L_b]\;\;| \;\forall a,b \in {\mathcal {U}} \},\; \;\;
\bar{{\mathcal {U}}}=\{\bar {A},A_a\;|\; 
\forall a \in {\mathcal {U}} \},\;(0.2)
$$
where $\{...\}$ is the linear span of elements in the braces.

Consider a direct sum
$${\mathcal {U}}\oplus {\mathcal {S}}\oplus\bar{{\mathcal {U}}}.\;\;(0.3)$$

Let $$a,b \in{\mathcal {U}} , \;\;\;\;\;S,S_1,S_2\in{\mathcal {S}} .\;\;\;\;\;\;(0.4)$$

Define
$$[a,b]=0,\;[\bar {A},A_b]=0,\;[A_a,A_b]=0,\;[S,a]=S(a),\;[S_1,S_2]=S_1S_2-S_2S_1,\;\;
(0.5)$$
$$[\bar {A},L_a]=A_a,\;[A_a,L_b]=A_{a*b}.\;\;\;(3.6)$$
\vspace {3 mm}
{\bf Theorem 1.} {\it The space (0.3) with the commutation relations (0.5) and (0.6)  is a 
graded Lie algebra
$$L(A)= U_{-1}\oplus U_{0}\oplus U_{1},\;\;\;\;(0.7)$$
where $U_{-1}={\mathcal {U}},\;U_0={\mathcal {S}},\;U_1=\bar{{\mathcal {U}}}.$}
\vspace {3 mm}

{\bf Remark  1.} The multiplication in the Jordan algebra $A$ can be
restored as
a double commutator $[[\bar {A},x],y]\;\;\;\forall x,y \in U_{-1}.$
\vspace {2 mm}

{\bf Example 1.} Let $A_n$ be a Jordan algebra of matrices of order $n$ with 
 operation $B*C=BC+CB$.

Then the Lie algebra $L(A_n)=A_{2n-1}$, i.e. linear Lie algebra of matrices of
order $2n$
$$A_{2n-1}=\left(
\begin{array}{c|c}
A&B\\
\hline
C&D\\
\end{array}\right),$$
(where $A,B,C,D$ are square matrices of order $n$) with the following grading
$$U_{-1}=\left\{ \left(
\begin{array}{c|c}
0&B\\
\hline
0&0\\
\end{array}\right)\right \},\;U_0=\left\{\left(
\begin{array}{c|c}
A&0\\
\hline
0&D\\
\end{array}\right)\right \},\;U_1=\left\{\left(
\begin{array}{c|c}
0&0\\
\hline
C&0\\
\end{array}\right)\right \}.$$


The matrix $\left(
\begin{array}{c|c}
  0&0\\
\hline
-E&0\\
\end{array}\right)$ 
plays the role of element $\bar A$ and
 the multiplication in the Jordan algebra $A_n$ can be restored as a double commutator
$$
\left[\left[\left(
\begin{array}{c|c}
  0&0\\
\hline
-E&0\\
\end{array}\right),\left(
\begin{array}{c|c}
0&x\\
\hline
0&0\\
\end{array}\right)\right],
\left(
\begin{array}{c|c}
0&y\\
\hline
0&0\\
\end{array}\right)\right]=
\left(
\begin{array}{c|c}
0&xy+yx\\
\hline
0&0\\
\end{array}\right).
$$


\begin{center} 

{\bf \S 2\  Jacobson's definition of a representation of a Jordan algebra
}

\end{center} 

Let $A$ be a Jordan algebra with multiplication $(x,y)\rightarrow xy$ and  $\phi : A
\rightarrow {\textstyle End\;}V$ be a linear map.

Consider a new algebra $\hat A$ on the space $A\oplus V$ with multiplication:
$$x*y=xy, \;u*v=0,\;x*v=v*x=\phi (x)v \;\;\;\forall ( x,y\in A, \;u,v \in V).\;\;\; (1.1)$$
In other words, $\hat A$ is a semi-direct sum of the given Jordan algebra  
 $A$ and an ideal $V$ with zero product.

{\bf Definition 1.} A linear map $\phi$ of a Jordan algebra $A$
into ${\textstyle End\;}V$ is called {\it a Jacobson 
representation } in $V$ if  $\hat A$ is a Jordan algebra. 

A similar definition in the Lie algebra case is equivalent to the ordinary one.
In a Jordan algebra case the notions of usual and Jacobson representations
are different.

The notions of  irreducibility and of equivalence of Jacobson representations can be given
in a natural way.

Let $A$ be a Jordan algebra with multiplication $x\ast y$. Let
$L_a(x)=a\ast x$. It is 
easy to show that the mapping $a \mapsto L_a$ is a Jacobson representation in the space of 
the algebra $A$.
It is natural to call this representation {\it the adjoint representation}.
Thus any Jordan algebra has at least one Jacobson representation.

\vspace {2 mm}

{\bf Theorem 2.} (N.Jacobson [2])
{\it There is only a finite  number of classes of irreducible Jacobson
  representations up to 
equivalence.} 
{\begin{center} 
{
\bf \S 3\  Generalized representations of a Jordan algebra
}\end{center} 

Let  $V$ be a linear space and  $a \in {\textstyle End\;}V$. We denote by  $*_a$ 
the multiplication on  $ {\textstyle End\;}V$ given by the  formula
 $$x*_ay=[[x,a],y]=xay+yax -axy -yxa. \;\;\;\; \;\;\;\;\ (2.1)$$

Let us call the space $ {\textstyle End\;}V$ with the multiplication $*_a$ {\it the algebra} 
$ {\textstyle End}_a\;V$.
 
{\bf Definition 2.} A homomorphism $\pi$ of an algebra $A$ with unit
$e$ into
some algebra 
${\textstyle End\;}_aV$
is called  {\it a generalized representation } of  $A$ if 
$$[a,[\pi (e),a]]=a.  \;\;\;(2.2)$$

{\bf Example 2.}
Let $A$ be a Jordan algebra and  $\phi$ be its usual representation in an $n$-dimensional 
linear space  $V$, i.e.,   let $\phi$ be a linear map of $A$ in  ${\textstyle End\;}V$ such that

$$\phi(x\cdot y)=\phi(x)\phi(y)+\phi(y)\phi(x). \; \;\;\;\;(2.3)$$

Consider a $2n$-dimensional space $\hat{V}=V_- \oplus V_+$, where $V_-$ and $V_+$ are 
copies of $V$. Consider in this space the following linear operators with matrices:
$$a=\left(
\begin{array}{c|c}
0&0\\
\hline
E&0\\
\end{array}\right),\;\;\;\;\pi(x)=\left(
\begin{array}{c|c}
0&\phi(x)\\
\hline
0&0\\
\end{array}\right).$$

Then it is easy to check that $\pi$ is a generalized representation of
the algebra $A$ into ${\textstyle End}_a\hat{V} $. Indeed,
$$\pi(x\cdot y)=\left(
\begin{array}{c|c}
0&\phi(x\cdot y)\\
\hline
0&0\\
\end{array}\right)=\left(
\begin{array}{c|c}
0&\phi(x)\phi(y)+\phi(y)\phi(x)\\
\hline
0&0\\
\end{array}\right)=$$
$$=[[\pi(x), a],\pi(y)]=\pi(x)*_a\pi(y),$$
and it is easy to see that  (2.2) is also fulfilled.

Thus to any usual representation $\phi$ of a Jordan algebra $A$ in a
linear space $V$
there corresponds a   generalized representation of this 
 algebra $A$ in the double
linear space $\hat{V}$.  We call the latter representation 
{\it a generalized representation of
 $A$ associated with a usual representation  $\phi$ }.

{\bf Example 3.}
Let $A$ be  a Jordan algebra  with a  multiplication $(x,y)\mapsto xy$
and  $\phi$ be
its Jacobson 
representation   in the space $V$. In other words, the algebra $\hat A$ defined on the 
space $A\oplus V$ with the multiplication (1.1) is a Jordan algebra.

Let us consider the Lie algebra $L(\hat A)= W_{-1}\oplus W_0\oplus W_1$. It is easy to see 
that $L(\hat A)=L(A)\oplus \hat V$, where $L(A)=U_{-1}\oplus U_0\oplus U_1$ and $\hat V =
V_{-1}\oplus V_0\oplus V_1$ is a graded commutative ideal in $L(\hat A)$ with $V_{-1}=V$. 

Denote by $a$ the element $x\mapsto [\bar A,x]\in {\textstyle
 End}(\hat V)$
where $\bar A$ is the element of $U_1\subset W_1$
such that $[[\bar A,x],y]=x*y\;\;\forall x,y \in   U_{-1} $, where $*$
 is the multiplication in the algebra $A$.
Then the mapping $\pi$ of $A=U_{-1}$ into $End(\hat V)$ defined by $\pi :b\mapsto [b,x]$, $x\in 
\hat V$ is a generalized representation of $A$ into algebra
${\textstyle End}_a(\hat V)$.

Thus to any Jacobson representation $\phi$ of a Jordan algebra $A$ in a
linear space $V$
there corresponds a   generalized representation of this 
 algebra $A$ in a
linear space $\hat{V}$, which we call {\it a generalized representation of
 $A$ associated with a Jacobson representation  $\phi$ }.

Let $\pi :A\to  {\textstyle End}_a\;V$ be a generalized representation of a Jordan algebra $A$ and 
 $W\subset V$ be an invariant subspace  with respect to $a$ and $\pi (A)$.
 Denote by $\phi : \pi(A)\cup\{a\}\to  {\textstyle End}\; W$ a homomorphism defined by  the 
restriction of elements of $\pi(A)\cup\{a\}$ to the space $W$.

Then the mapping $\pi_W :A\to  {\textstyle End}_{\phi (a)} W$ defined by $\pi_W(x)=\phi(\pi(x))$ is a 
generalized representation of $A$ called {\it a subrepresentation  } of $\pi$.

{\bf Definition 3.} A generalized representation $\pi$ is called {\it irreducible} if $\pi$
has no subrepresentations.

{\bf Definition 4.} Two  generalized representations $\pi:A\to
{\textstyle  End}_a\; V$ and $\nu 
:A\to {\textstyle  End}_b\; W$  are called {\it equivalent} if there exists an isomorphism $f:V\to W$ 
such that
$f\pi (x)=\nu f(x)$ if $x\in V$ and $fa=bf$.

Let $A$ be a Jordan algebra with unit and $\pi$ be a representation of Lie algebra $L(A)$. 
Denote by $\Phi (\pi)$ the restriction of $\pi$ to the subspace $U_{-1}$.

{\bf Lemma.} {\it  Let $\pi$ be a representation of  $L(A)$ in 
the space $V$.  The mapping $\Phi $ associates to it a generalized representation of Jordan algebra $A$
in $ {\textstyle End}_{\pi (\bar A)}\;V$.}

The proof  is very simple:
$$\Phi (\pi)(A(x,y))={\pi}([[[\bar {A},x],y])=[[{\pi}([\bar 
{A}),{\pi}(x)],{\pi}(y)]
=\Phi (\pi(x))*_a\Phi (\pi(y)),\;(2.4)$$
where $a= {\pi}(\bar {A}).$
The second condition (2.2) is also fulfilled:
$$[ {\pi}(\bar {A}),[\pi (e), {\pi}(\bar {A})]]=
{\pi}([\bar {A},[e,\bar {A}]]={\pi}(\bar {A}),  \;\;\;(2.5)$$
 because $[\bar {A},[e,\bar {A}]]=\bar {A}$  in the algebra Lie $L(A)$.

 The following statement allows to describe all irreducible
 representations of a Jordan 
algebra with unit.
\vspace {3 mm}

{\bf Theorem 3.} {\it  The mapping  $ \Phi $
  defines a one-to-one correspondence between 
 irreducible linear representations    of Lie algebra  $L(A)$ 
and irreducible generalized representations of Jordan algebra $A$.}

\begin{center} 

{\bf \S 4\
The  graded representations of graded Lie algebras \\and the order of
generalized representations }
\end{center} 

Let $G$ be a graded Lie algebra
$$G= U_{-k}\oplus \cdots  \oplus U_{-1}\oplus U_{0} \oplus U_{1}\oplus
\cdots  \oplus U_{m},\;\;\;\;(3.1)$$

Let $\pi$ be a representation of $G$  in the linear space $V$. {\it An $l$-grading of 
$\pi$} is 
a presentation of $V$ as a direct sum:
$$V=  V_{1} \oplus V_{2}\oplus \cdots  \oplus V_{l},\;\;\;\;(3.2)$$
such that
 $$\pi( U_{i})V_{j}\subset V_{i+j}.\;\;\;\;(3.3)$$

A representation $\pi$ of $G$  in the linear space $V$ equipped with
an $l$-grading is called an
{\it $l$-graded representation}.
The number $l$ is called {\it the length} of the graded representation $\pi$.

 {\bf Example 4.} Let $A_n$ be a Jordan algebra of matrices of order $n$ with 
the operation $B*C=BC+CB$.
 Consider the  graded Lie algebra $A_{2n-1}=U_{-1}\oplus U_{0} \oplus U_{1}$
 from example~1.
The  representation of this algebra given  there is 2-graded:

$$V=  V_{1} \oplus V_{2},$$ 
where $ V_{1}$ is the linear subspace spanned by the first $n$
coordinate vectors and   $V_{2} $ is the linear subspace spanned by the
last $n$ coordinate vectors. 

{\bf Definition 5.} A generalized representation of a Jordan algebra $A$
 is {\it  of order $l$}
 if it  has the form $\Phi (\pi )$, where $\pi$  is an $l$-graded
 representation of the Lie algebra $L(A)$.

{\bf  Remark 2.} As Examples 2 and~3 show, a generalized representation of 
Jordan algebra $A$  associated with the usual (the Jacobson)
representation is of order~2 (order~3).
\vspace {3 mm}

{\bf  Theorem 4.}  {\it To find 
   all usual (all Jacobson) representations  of a Jordan
   algebra $A$ it is enough to find all $2$-graded ($3$-graded) 
  representations   $ {\pi}$
of the Lie algebra  $L(A)=U_{-1}\oplus U_{0} \oplus U_{1}$  and consider
$\Phi (\pi )$.

        In particular, a Jordan   algebra $A$ is special iff  the Lie algebra
  $L(A)$ has a $2$-graded representation.}
\vspace {3 mm}

{\bf Corollary } (See [1],[2],[3]).
 {\it The exceptional Jordan algebra $E_3$ has no usual representations.}

This classical result follows also from the following simple fact:
the graded Lie algebra $E_7=L(E_3)=U_{-1}\oplus U_{0} \oplus U_{1}$
has no 2-graded representations.

\begin{center} 
{\bf \S 5\  Classification of irreducible generalized representations 
\\of  simple Jordan algebras.}
\end{center} 

To classify irreducible generalized representations of  simple Jordan
algebras we first will 
describe all gradings of  irreducible representations of graded simple Lie algebras.

Recall that a grading of a semisimple Lie algebra $G$ is defined by a linear
function $f(x)$ on the dual space $H^*$ of a Cartan subalgebra
$H\subset G$ taking nonnegative integer 
values at simple roots of the algebra  $G$ (see [6],[7]).

{\bf Theorem 5.}
 {\it Let $G=\sum\limits _{i=-k} ^{i=k} U_{i}$ be a semisimple Lie algebra
graded by a function $f$. Then a linear space $V$
of a   representation  $\pi$ of $G$ with the highest weight $\Lambda$ can be 
equipped with a grading such that  $\pi$ becomes an  $l$-graded
representation with $$l=f(\Lambda)+f(i(\Lambda))+1,$$
where $i$ is a Tits involution on simple roots of~$G$.
The length $l$ is uniquely defined.
}
\vspace {3 mm}

{\bf Remark 3.} Recall that the Tits involution is defined by a nontrivial involutive 
automorphism of the Dynkin diagram of $G$ in the cases $G=A_n,D_{2n},E_6$ and the identity
in other cases.

  We will apply Theorem~5 to the  special cases when
  $G=L(A)= U_{-1}\oplus U_{0} \oplus 
U_{1}$. In all these cases the function $f$ is equal to~$1$ at some
simple root
$\alpha _j$  
and is equal  to zero at other simple roots, i.e. $f(\Lambda )
=\lambda _j$,
 where $\Lambda =\sum\limits_1^n\lambda _i \alpha _i$ is the
 decomposition of
 $\Lambda$
into a linear combination of simple roots $ \alpha _i $.

Let $A$ be a simple Jordan algebra and $\pi _\Lambda $
be an irreducible representation of the Lie algebra $L(A )$
with the highest weight $\Lambda$.

The following theorems follow from Theorems 4 and~5.
\vspace {2 mm}
{\bf Theorem 6.} {\it The finite dimensional irreducible generalized representations of a 
simple Jordan algebra $A$ over a field $C$ have the form $\Phi (\pi _\Lambda )$.

The orders $l$ of these  representations  are given by the following list. 

1) $A=A_n$, $L(A)=A_{2n-1}$, $f(\Lambda )=\lambda_n$, $l=2\lambda_n+1$.

2) $A=B_n$, $L(A)=C_n$, $f(\Lambda )=\lambda_n$, $l=2\lambda_n+1$.

3)  $A=C_n$, $L(A)=D_{2n}$, $f(\Lambda )=\lambda_{2n}$, $l=\lambda_{2n}+\lambda_{2n-1}+1$.

4) $A=D_{2n}$, $L(A)=D_n$, $f(\Lambda )=\lambda_1$, $l=2\lambda_1+1$.

5) $A=D_{2n+1}$, $L(A)=B_n$, $f(\Lambda )=\lambda_1$, $l=2\lambda_1+1$.

6) $A=E_3$, $L(A)=E_7$, $f(\Lambda )=\lambda_7$, $l=2\lambda_7+1$.
}

\vspace {3 mm}

{\bf Theorem 7.}  {\it Let $A$ be a simple Jordan algebra.
  There is a finite number  of irreducible generalized representations
  of a given order~$l$.

  In particular, there is a finite number of usual representations ($l=2$)
and Jacobson representations ($l=3$).}

\begin{center} 
{\bf \S 6\ An interpretation of the generalized representations} 

\end{center}

First we will define a  series of algebras $A_{k_1k_2\cdots
  k_{l}}$. Let $U_{st}$ be the space of $s \times t$
matrices.

The space of algebra  $A_{k_1k_2\cdots
  k_{l}}$ is the direct sum $$A_{k_1k_2}\oplus A_{k_2k_3}\oplus
\cdots
\oplus A_{k_{l-1}k_{l}},$$
i.e. the elements of the algebra are the tuples
$$x=(x_1,x_2,\cdots ,  x_{l-1}),\;\;\;\;\;x_i \in  U_{k_{i}k_{i+1}}.$$

Let us fix $a=(a_1,a_2,\cdots ,  a_{l-1}).$ The multiplication  $*_a$
in the algebra  $A_{k_1k_2\cdots
  k_{l}}$ is given by the formula
$$(x*_a y)_i=x_ia^t_i y_i +y_ia^t_ix_i-a^t_{i-1}x_{i-1} y_i -y_i x_{i+1}a^t_{i+1},$$
where $a^t$ is the transpose of matrix $a$.

{\bf Remark 4.} The algebra  $A_{kk}$ is the Jordan algebra of matrices
of order $k$. Indeed
$$x*_a y=xa^ty +ya^tx.$$
It is easy to see that this multiplication is isomorphic to the 
multiplication $xy +yx$
 if $a$ is not a degenerate matrix.
\vspace {2 mm}

{\bf Theorem 8.} {\it There is a one-to-one correspondence between generalized 
representations $\pi$ of order $l$ of Jordan algebra $A$
and homomorphic embeddings   $ {\pi}^*$ of 
 $A$ in algebras $A_{k_1k_2\cdots
  k_{l}}$.}

\begin{center}
\bf References
\end{center}
\vskip 3mm
\rm
\begin{enumerate}

\item  
 P. Jordan, J. von Neumann, E. Wigner, 
On the algebraic generalization of the quantum mechanical formalism,
Ann. of Math. (2)  35 (1934), 29-64.

\item
A.Albert,
on a certain algebra of quantum mechanics, Ann. of Math. (2)  35 (1934), 65-73.

\item  
N. Jacobson,
Structure and representations of Jordan algebras,
Amer.Math. Soc. Vol. 39 (1968), 149-170.

\item
I. L. Kantor,
Classification of irreducible transitive differential groups, Dokl. Akad.
 Nauk SSR,
        Vol. 158, No. 6, 1271-1274 (1964). 

 \item
I. L. Kantor, Transitive differential groups and invariant connections
on homogeneous spaces.
  Trudy Sem. Vector Tensor Anal., Vol. X111, 227-266
(1966).

\item
V.G. Kac, Classification of simple Z-graded Lie superalgebras and
simple Jordan superalgebras. Comm. in algebra, 5(13),
(1977), 1375-1400.

\item
I. L. Kantor,
Some generalizations of Jordan 
algebras,
Trudy Sem.Vektor.Tensor.Anal.16 (1972), 407-499 (Russian).

\end{enumerate}
\end{document}